\documentclass[aop,MSNbibl,citesort,dvips]{arximspdf}
\usepackage{graphicx}
\doi{10.1214/10-AOP634}
\volume{39}
\issue{5}
\pubyear{2011}
\firstpage{1815}
\lastpage{1843}

\makeatletter

\let\sv@tabnotetext\tabnotetext
  \let\sv@tabnotemark@fmt\tabnotemark@fmt
   \long\def\legend#1{{\let\tabnote@indent\leavevmode\sv@tabnotetext[]{}{#1}}}

\def\bptnote#1{}

\newcommand{\underset}[2]{\mathop{#1}_{#2}}

\newtheorem{theorem}{Theorem}
\newtheorem{prop}[theorem]{Proposition}
\newtheorem{lemma}[theorem]{Lemma}
\newproclaim{defi}[theorem]{Definition}
\newtheorem{cor}[theorem]{Corollary}
\newproclaim{rmk}[theorem]{Remark}

\newcommand{\eps}{\varepsilon}

\newcommand{\var}{\operatorname{var}}
\newcommand{\cov}{\operatorname{cov}}

\newcommand{\E}{\mathbb{E}}
\newcommand{\PP}{\mathbb{P}}

\newcommand{\ggnew}{|\Gamma|}
\newcommand{\g}{\Gamma}

\newcommand{\cA}{\mathcal A}
\newcommand{\cF}{\mathcal F}
\newcommand{\cS}{\mathcal S}
\newcommand{\Poi}{\operatorname{Poisson}}
\newcommand{\tm}{t_{\mathrm{mix}}}

\makeatother

\begin{document}
\begin{frontmatter}

\title{Mixing times for random $\bolds k$-cycles and coalescence-fragmentation chains}
\runtitle{Mixing times for random $k$-cycless}

\begin{aug}
\author[A]{\fnms{Nathana\"{e}l} \snm{Berestycki}\thanksref{t1}\ead[label=e1]{N.Berestycki@statslab.cam.ac.uk}%
\ead[label=u1,url]{http://www.statslab.cam.ac.uk/\textasciitilde beresty/}},
\author{\fnms{Oded} \snm{Schramm}} and
\author[C]{\fnms{Ofer} \snm{Zeitouni}\corref{}\thanksref{t3}\ead[label=e3]{zeitouni@math.umn.edu}%
\ead[label=u3,url]{http://www.wisdom.weizmann.ac.il/\textasciitilde zeitouni/}}

\runauthor{N. Berestycki, O. Schramm and O. Zeitouni}
\affiliation{Cambridge University, and Weizmann Institute and
University of Minnesota}
\dedicated{Dedicated to the memory of Oded Schramm}
\address[A]{N. Berestycki\\
CMS---Wilberforce Rd.\\
Cambridge University\\
Cambridge CB3 0WB\\
United Kingdom\\
\printead{e1}\\
\printead{u1}} 
\address[C]{O. Zeitouni\\
Faculty of Mathematics\\
Weizmann Institute\\
POB 26, Rehovot 76100\\
Israel\\
and\\
School of Mathematics\\
University of Minnesota\\
206 Church St. SE\\
Minneapolis, Minnesota 55455\\
USA\\
\printead{e3}\\
\printead{u3}}
\end{aug}

\thankstext{t1}{Supported in part by EPSRC Grant EP/GO55068/1.}
\thankstext{t3}{Supported in part by NSF Grant DMS-08-04133 and by a
grant from the Israel
Science Foundation.}

\received{\smonth{11} \syear{2010}}

%
\begin{abstract}
Let $\mathcal{S}_n$ be the permutation group on $n$ elements,
and consider a~random walk on $\mathcal{S}_n$
whose step distribution is uniform on $k$-cycles.
We prove a well-known
conjecture that the mixing time of this process is
$(1/k) n \log n$, with threshold of width linear in $n$.
Our proofs are elementary and purely probabilistic,
and do not appeal to the representation theory of $\mathcal{S}_n$.
\end{abstract}

%
\begin{keyword}[class=AMS]
\kwd{60B15}
\kwd{60J27}.
\end{keyword}
\begin{keyword}
\kwd{Mixing times}
\kwd{coalescence}
\kwd{cutoff phenomena}
\kwd{random cycles}
\kwd{random transpositions}.
\end{keyword}

\end{frontmatter}

\section{Introduction}

\subsection{Main result}

Let $\mathcal{S}_n$ be the group of permutations of $\{1, \ldots, n\}$.
Any permutation $\sigma\in\mathcal{S}_n$ has a unique cycle
decomposition, which partitions the set $\{1,\ldots, n\}$ into orbits
under the natural action of $\sigma$. The cycle structure of $\sigma$
is the integer partition of $n$ associated with this set partition, in
other words, the ordered sizes of the cycles (blocks of the partition)
ranked in decreasing size. It is customary not to include the fixed
points of $\sigma$ in this structure. For instance, the permutation
\[
\sigma= \pmatrix{
1 & 2 & 3 & 4 & 5 & 6 & 7 \cr
4 & 2 & 6 & 7 & 3 & 5 & 1
}
\]
has 3 cycles, $(1 \enskip4 \enskip7) (2) (3 \enskip6 \enskip5)$, so
its cycle structure is
$(3,3)$ (and one fixed point which does not appear in this structure).
A conjugacy class $\Gamma\subset\mathcal{S}_n$~is~the set of permutations
having a given cycle structure.\vadjust{\eject} Let $|\Gamma|$ denote the support of
$\Gamma$,
that is, the number of nonfixed-points of any permutation $\sigma\in
\Gamma$.
In what follows we deal with the case where $\Gamma$
consists of a single $k$-cycle,
in which case $|\Gamma|=k$ (see, however,
Remark \ref{rem-2}).
It is well known and easy to see that in this case,
if $k$ is even, then $\Gamma$ generates $\mathcal{S}_n$, while if $k
> 2$
is odd, then $\Gamma$ generates the alternate group $\mathcal{A}_n$ of
even permutations. Let $(\pi_t, t\ge0)$ be the continuous-time random
walk associated with $(\cS_n, \Gamma)$. That is, let $\gamma_1,
\gamma
_2, \ldots$ be a sequence of i.i.d. elements uniformly distributed on
$\Gamma$, and let $(N_t, t \ge0)$ be an independent Poisson process
with rate 1; then we take
%
\begin{equation}
\pi_t = \gamma_1 \circ\cdots\circ\gamma_{N_t},
\end{equation}
where $\gamma\circ\gamma'$ indicates the composition of the
permutations $\gamma$ and $\gamma'$. $(\pi_t , t\ge0)$ is a Markov
chain on $\cS_n$ which converges to the uniform distribution $\mu$
on~$\cS_n$ when $\ggnew$ is even, and to the uniform distribution on
$\cA
_n$ when $\ggnew>2$ is odd. In any case we shall write $\mu$ for that
limiting distribution. We shall be interested in the mixing properties
of this process as $n\to\infty$, as measured in terms of the total
variation distance. Let $p_t(\cdot)$ be the distribution of $\pi_t$ on
$\cS_n$, and let $\mu$ be the invariant distribution of the chain. Let
\[
d(t) = \|p_t(\cdot) - \mu\| = \frac12\sum_{\sigma\in\cS_n}
|p_t(\sigma
)- \mu(\sigma)|,
\]
where $d(t)$ is the total variation distance between the state of the chain
at time $t$ and its limiting distribution $\mu$.
(Below, we will also use the notation $\|X-Y\|$ where $X$ and $Y$ are
collections of random variables with laws $p_X,p_Y$ to mean $\|p_X-p_Y\|$.)

The main goal of this paper is to prove that the chain exhibits a sharp
cutoff, in the sense that $d(t)$ drops abruptly from its maximal value
1 to its minimal value 0 around a certain time $\tm$, called the mixing
time of the chain. (See \cite{diaconis} or \cite{LPW} for a general
introduction to mixing times.) Note that if $\Gamma$ is a fixed
conjugacy class of $\cS_n$ and $m>n$, $\g$ can also be considered a
conjugacy class of $\cS_m$ by simply adding $m-n$ fixed points to any
permutation $\sigma\in\g$. With this in mind, our theorem states the
following:
\begin{theorem} \label{T:mix}
Let $k\geq2$ be an integer, and let $\Gamma_k$ be the
conjugacy class
of $\cS_n$ corresponding to $k$-cycles. The continuous time
random walk $(\pi_t, t\geq0)$
associated with
$(\cS_n,\Gamma_k)$ has a cutoff at time
$\tm:= (1/k) n \log n$, in the sense that for any $\eps>0$,
there exist $N_{\eps,k},C_{\eps,k}>0$
large enough so that for all $n \ge N_{\eps,k}$,
%
\begin{eqnarray}
\label{Tlb}
d(\tm-C_{\eps,k}n) &>& 1- \eps, \\
\label{Tub}
d(\tm+C_{\eps, k}n ) &<& \eps.
\end{eqnarray}
\end{theorem}

As explained in Section
\ref{subsec-background}
below, this result solves a well-known conjecture formulated by several
people over the course of the years.\vadjust{\eject}
\begin{rmk}
\label{rem-2}
Theorem \ref{T:mix} can be extended, without a significant change in
the proofs, to cover the case of general fixed conjugacy classes
$\Gamma$, with $k=|\Gamma|>2$ independent of $n$. In order to alleviate
notation, we present here only the proof for $k$-cycles.
A more
delicate question, that we do not investigate, is what growth
of $k=k(n)$ is allowed so that Theorem \ref{T:mix} would still be true
in the form
%
\begin{eqnarray}
\label{Tlb1}
d\bigl(\tm(1-\delta)\bigr) &>& 1- \eps, \\
\label{Tub1}
d\bigl(\tm(1+\delta)\bigr) &<& \eps?
\end{eqnarray}
The lower bound in (\ref{Tlb1}) is easy. For the upper bound in
(\ref{Tub1}), due to the birthday problem, the case $k= o(\sqrt{n})$
should be fairly similar to the arguments we develop below, with
adaptations in several places, for example, in the argument
following~(\ref{eq-070110}); we have not checked the details. Things
are likely to become more delicate when $k $ is of order $\sqrt{n}$ or
larger. Yet, we conjecture that (\ref{Tub1}) holds as long as
\mbox{$k=o(n)$}.
\end{rmk}

\subsection{Background}
\label{subsec-background}

This problem has a rather long history, which we now sketch. Mixing
times of Markov chains were studied
independently by Aldous \cite{aldous-mix} and by Diaconis and
Shahshahani \cite{d-sh} at around the same time, in the early 1980s.
Diaconis and Shahshahani \cite{d-sh}, in particular, establish the
existence of what has become known as the \textit{cutoff phenomenon} for
the composition of random transpositions. Random transpositions is
perhaps the simplest example of a random walk on $\cS_n$ and is a
particular case of the walks covered in this paper, arising when the
conjugacy class $\g$ contains exactly all transpositions.
The authors of \cite{d-sh}
obtained a version of
Theorem \ref{T:mix} for this particular case (with explicit choices of
$C_{2,\eps}$ for a given $\eps$).
As is the case here, the hard part of the result is the upper-bound
(\ref{Tub}). Remarkably, their solution involved a connection with the
representation theory of~$\cS_n$, and uses rather delicate estimates on
so-called character ratios.

Soon afterwards, a flurry of
papers tried to generalize the results of \cite{d-sh} in the direction
we are taking in this paper, that is, when the step distribution is
uniform over a fixed conjugacy class $\g$. However, the
estimates on character ratios that are needed become harder and
harder as $\ggnew$ increases. Flatto, Odlyzko and Wales \cite{FOW},
building on earlier work of Vershik and Kerov \cite{vk},
obtained finer estimates on character ratios and were able to
show that mixing must occur before $(1/2) n \log n$ for $\ggnew$ fixed,
thus giving another proof of the Diaconis--Shahshahani result when
$\ggnew=2$.
(Although this does not appear explicitly in \cite{FOW},
it is recounted in Diaconis's book \cite{diaconis}, page 44.)
Improving further the estimates on character ratios, Roichman
\cite{roichman1,roichman2} was able to prove a weak version of Theorem \ref
{T:mix}, where it is shown that $d(t)$ is small if
$t> C\tm$ for some large enough $C>0$. In his result, $\ggnew$ is
allowed to grow to infinity as fast as $(1-\delta) n$ for any $\delta>0$.
To our knowledge,
it is in \cite{roichman1} that Theorem \ref{T:mix} first formally
appears as a conjecture, although we have no doubt that it had been
privately made before. (The lower bound for random transpositions,
which is based on counting the number of fixed points in $\pi_t$, works
equally well in this context and provides the
conjectured correct
answer in all cases.)
Lulov \cite{lulov} dedicated his Ph.D. thesis to the problem, and
Lulov and Pak \cite{lulov-pak} obtained a partial proof of the
conjecture of Roichman, in the case where $\ggnew$ is very large, that
is, greater than $n/2$. More recently, Roussel \cite{roussel1} and
\cite{roussel2} made some progress in the small $\ggnew$ case, working out
the character ratios estimates to treat the case where $\ggnew\le6$.
Saloff-Coste, in his survey article (\cite{lsc}, Section~9.3) discusses
the sort of difficulties that arise in these computations and states
the conjecture again. A summary of the results discussed above is also
given. See also
\cite{SZ}, page 381, where work in progress of Schlage-Puchta that
overlaps the result in Theorem \ref{T:mix} is mentioned.

\subsection{Structure of the proof}

To prove Theorem \ref{T:mix}, it suffices to look at the cycle
structure of $\pi_t$ and check that if $N_t(i)$ is the number of cycles
of $\pi_t$ of size $i$ for every $i\ge1$, and if $t\ge\tm+
C_{k,\eps
}n$ then the total variation distance between $(N_t(i))_{1\le i \le n}$
and $(N(i))_{1\le i \le n}$ is close to 0, where $(N(i))_{1\le i \le
n}$ is the cycle distribution of a random permutation sampled from $\mu
$. We thus study the dynamics of the cycle distribution of $\pi_t$,
which we view as a certain coagulation--fragmentation chain. Using ideas
from Schramm \cite{schramm}, it can be shown that large cycles are at
equilibrium
much before $\tm$, that is, at a time of order $O(n)$. Very informally
speaking, the idea of the proof is the following. We focus for a moment
on the case $k=2$ of random transpositions, which is the easiest to
explain. The process $(\pi_t,t\ge0)$ may be compared to an Erd\H
{o}s--R\'{e}nyi random graph process $(G_t,t\ge0)$ where random edges are
added to the graph at rate 1, in such a way that the cycles of the
permutation are subsets of the connected components of $G_t$. Schramm's
result from \cite{schramm} then says that, if $t=cn$ with $c>1/2$ (so that
$G_t$ has a giant component), then the macroscopic cycles within the
giant component have relaxed to
equilibrium.
By an old result of Erd\H{o}s and R\'{e}nyi, it takes time $t=\tm+
C_{k,\eps}n$ for $G_t$ to be connected with probability greater than
$1-\eps$. By this point the giant component encompasses every vertex
and thus, extrapolating Schramm's result to this time, the macroscopic
cycles of $\pi_t$ have the correct distribution at this point. A
separate and somewhat more technical argument is needed to deal with
small cycles.

More formally, the proof of Theorem \ref{T:mix} thus proceeds in two
main steps. In the first step,
presented in Section \ref{smallcycles} and culminating
in Proposition \ref{prop-small},
we show that after time
$\tm+c_{\eps,k}n$, the distribution of \textit{small cycles}
is close (in variation distance) to the invariant measure, where a
\textit{small cycle} means that it is smaller than a suitably chosen threshold
approximately equal to $n^{7/8}$. This is achieved
by combining a queueing-system argument (whereby initial discrepancies
are cleared by time slightly larger than
$\tm$ and equilibrium is achieved) with a priori
rough estimates on the decay of mass in small cycles (Section
\ref{subsec-verif}). In the second step, contained in
Section \ref{sec-schramm},
a variant of Schramm's coupling
from \cite{schramm} is presented, which allows us to
couple the chain after time
$\tm+c_{\eps,k}n$
to a chain started from equilibrium,
within time of order $n^{5/8}\log n$, if all small cycles agree initially.

\section{Small cycles}
\label{smallcycles}

In this section we prove the following proposition. Let $(N_i(t))_{1
\le i \le n}$ be the number of cycles of size $i$ of the permutation
$\pi_t$, where $(\pi_t,t\ge0)$ evolves according to random $k$-cycles
(where $k\ge2$), but does not necessarily start at the identity
permutation. Let $(Z_i)_{i=1}^n$ denote independent Poisson random
variables with mean $1/i$.

Fix $ 0 < \chi< 1 $ and let
$K=K(n)$ be the closest dyadic integer to
$n^\chi$.
We think of cycles smaller than $K$ as being small, and big otherwise.
Let
$I_j=\{i\in\mathbb{Z}\dvtx i\in[2^j,2^{j+1})\}$, $L_j=|I_j|=2^j$ and
%
\begin{equation}
\label{eq-mj}
M_j(t)=\sum_{i\in I_j} N_i(t) .
\end{equation}
Introduce the stopping time
%
\begin{equation}\label{D:tau}
\tau= \inf\{t\ge0\dvtx \exists0 \le j \le\log_2K+1,
M_j(t) > (\log n)^6/2\}.
\end{equation}
Therefore, prior to $\tau$,
the total number of small cycles in
each dyadic strip
$[2^j$, $2^{j+1})$ $(j\le1+\log_2 K )$ never exceeds $(\log n)^6/2$.
\begin{prop}
\label{P:step2}
Suppose that
%
\begin{equation}\label{assumptionstep2}
\PP(\tau< n \log n) \longrightarrow0
\end{equation}
as $n \to\infty$,
and that initially,
%
\begin{equation}
\label{assumptionstep2prime}
M_j(0) \le D \log(j+2)
\end{equation}
for all $0\le j\le\log_2 \log n$, for some $D>0$ independent of $j$
or $n$.
Then for any sequence $t=t(n)$ such that $t(n)/n \to\infty$ as $n\to
\infty$ and $t(n) \le n \log n$,
%
\[
\| (N_i(t))_{i=1}^{K} - (Z_i)_{i=1}^{K}\|
\longrightarrow0.
\]
\end{prop}

In particular,
under the assumptions of Proposition \ref{P:step2}, for any $\eps>0$
there is a $c_{\eps,k}>0$ such that for all $n$ large,
\[
\| (N_i(c_{\eps,k}n))_{i=1}^{K} -
(Z_i)_{i=1}^{K}
\| <\eps.
\]
%

In Sections \ref{subsec-verif}
and \ref{subsec-conclusionrev},
Proposition \ref{P:step2} is
applied to the chain after time roughly $\tm=(n\log n)/k$, at which
point the initial conditions $M_j(0)$ satisfy~(\ref{assumptionstep2prime})
(with high probability).
\begin{pf*}{Proof of Proposition \ref{P:step2}}
The proof of this proposition relies on the
analysis of the dynamics of the small cycles, where each step
of the dynamics corresponds to an application of a $k$-cycle, by
viewing it as a coagulation--fragmentation process.
To start with, note that every $k$-cycle may decomposed as a product of
$k-1$ transpositions
\[
c=(x_k, \ldots, x_1) = (x_k, x_{k-1}) \cdots (x_2, x_1).
\]
Thus the application of a $k$-cycle may be decomposed into the
application of $k-1$ transpositions: namely, applying $c$ is the same
as first applying the transposition $(x_1, x_2)$ followed by
$(x_{2},x_{3})$ and so on until $(x_{k-1},x_k)$. Whenever one of those
transpositions is applied, say $(a, b)$, this can yield either a
fragmentation or a coagulation, depending on whether $a$ and $b$ are in
the same cycle or not at this time. If they are, say if $b= \sigma
^i(a)$ (where $i\ge1$ and~$\sigma$ denotes the permutation at this
time), then the cycle $C$ containing $a$ and~$b$ splits into
$(a,\ldots, \sigma^{i-1}(a))$ and everything else, that is, $(b, \ldots, \sigma
^{|C| - i}(b))$. If they are in different cycles $C$ and $C'$ then the
two cycles merge.

To track the evolution of cycles,
we color the cycles with different colors (blue, red or black)
according (roughly) to the following rules. The blue cycles will be the
large ones, and the small ones consist of red and black. Essentially,
red cycles are those which undergo a ``normal'' evolution, while the
black ones are those which have experienced some kind of error. By
``normal evolution,'' we mean the following: in a given step,
one
small cycle is generated by fragmentation of a blue cycle. It is the
first
small cycle that is involved in this step. In a later step of the
random walk, this cycle coagulates with a large cycle and thus becomes
large again. If at any point of this story, something unexpected
happens (e.g., this cycle gets fragmented instead of coagulating with a
large cycle, or coagulates with another small cycle) we will color it
black. In addition, we introduce ghost cycles to compensate for this
sort of error.

We now describe this procedure more precisely. We start by coloring
every cycle of the permutation $\sigma(t)$ which is larger than $ K$
blue. We denote by~$\theta(t)$ the fraction of mass contained in blue
cycles, that is,
%
\begin{equation}
\theta(t)= \frac1n \sum_{i = K+1 }^n i N_i(t).
\end{equation}
Note that by definition of $\tau$,
%
\begin{equation}
1- \frac{K}n (\log n)^6 \le\theta(t) \le1
\end{equation}
for all $t \le\tau$.

We now color the cycles which are smaller than $K$ either red or black
according to the following dynamics. Suppose we are applying a certain
$k$-cycle $c=(x_k, \ldots, x_1)$, which we write as a product of $k-1$
transpositions
%
\begin{equation}
c=(x_k, \ldots, x_1)=(x_k, x_{k-1}) \cdots(x_2,x_1)
\end{equation}
(note that we require that $x_i\neq x_{j}$ for $i\neq j$).\vspace*{8pt}

\textit{Red cycles.}
Assume that a blue cycle is fragmented and one of the pieces is small,
and that this transposition is the first one in
the application of the $k$-cycle $(x_1, \ldots, x_k)$ to
involve a small cycle.
In that case (and only in that case), we color it red.
Red cycles may depart through coagulation or fragmentation.
A coagulation with a blue cycle, if it is the first in the step and no
small cycles were created in this step prior to it, will be
called \textit{lawful}. Any other departure will be called \textit{unlawful}.
If a blue cycle breaks up in a~way that would create a red cycle
and both cycles created
are small (which may happen if the size of the cycle is between $K$ and $2K$),
then we color the smaller one red and the larger one black, with a
random rule in the case of ties.\vspace*{8pt}

\textit{Black cycles.}
Black cycles are created in one of two ways. First,
any red cycle that
departs in an unlawful fashion and stays small becomes black.
Further,
if the transposition $(a,b)$ is not the first transposition
in this step to create a small cycle from a blue cycle, or if
it is but a previous transposition in the step involved a small cycle,
then the small
cycle(s) created
is colored black. Now, assume that $(a, b)$ involves only cycles
which are smaller than $K$: this may be a fragmentation producing
two new cycles, or a merging of two cycles producing one new cycle.
In this case, we color the new cycle(s) black,
no matter what the initial color of the cycles,
except if this operation is a coagulation \textit{and} the size of
this new cycle exceeds $K$, in which case it is colored blue again.
Thus, black cycles are created through either coagulations
of small parts or fragmentation of either small or large parts,
but black cycles disappear only through coagulation.

We aim to analyze the dynamics of the red and black system,
and the idea is that the dynamics of this system are essentially
dominated by that of the red cycles, where the occurrence of black
cycles is an error that we aim to control.\vspace*{8pt}

\textit{Ghosts.} Let $R_i(t), B_i(t)$ be the number of red
and black cycles, respectively, of size $i$ at time $t$.
It will be helpful to introduce another type of cycle, called ghost
cycles, which are nonexisting cycles which we add for counting purposes:
the point is that we do not want to touch more than one red cycle in
any given step. Thus, for any red cycle departing in an unlawful way,
we compensate it by creating a ghost cycle of the same size. For
instance, suppose two red cycles $C_1$ and $C_2$ coagulate (this could
form a blue or a~black cycle). Then we leave in the place of $C_1$ and
$C_2$ two ghost cycles $C'_1$ and $C'_2$ of sizes identical to $C_1$
and $C_2$.

An exception to this rule is that if, during a step,
a transposition creates a small red cycle by
fragmentation of a blue cycle, and later within the same step
this red cycle either is immediately
fragmented again in the next transposition
or coagulates with another red or black cycle and remains small,
then it becomes black as above but we do not leave a ghost in its place.

Finally, we also declare that every ghost cycle of size $i$ is killed
independently of anything else at an instantaneous rate which is
precisely given by~$i \mu(t)$, where $\mu(t)$ is a random nonnegative
number (depending on the
state of the system at time $t$)
which will be defined below in (\ref{mu}) and corresponds to the rate
of lawful departures of red cycles.

To summarize, we begin at time $0$
with all large cycles colored blue and all small cycles
colored red. For every step consisting of $k$ transpositions,
we run the following algorithm
for the coloring of
small cycles and creation of ghost cycles (see Table \ref{table1}).

\begin{table}[t]
\caption{Coloring algorithm for small cycles, and creation of ghost cycles}
\label{table1}
\rule{\tablewidth}{0.5pt}\vspace*{-5pt}

\begin{itemize}
\item[$\bullet$ (I)] If the transposition is a fragmentation, go to (F);
otherwise, go to (C).
\item[$\bullet$ (F)] If the fragmentation is of a small cycle $c$
of length $\ell$, go to (FS); otherwise, go to (FL).
\item[$\bullet$ (FS)] Color the
resulting small cycles black. Create a ghost cycle of length
$\ell$, except if
$c$ was created in the previous transposition of
the current step and is red. \textit{Finish}.
\item[$\bullet$ (FL)] If the fragmentation
creates one or two
small cycles, and this transposition is the first in
the step to either create or involve a small cycle,
color the smallest small cycle created red. All other small
cycles created are colored black. Do not
create ghost cycles. \textit{Finish}.
\item[$\bullet$ (C)] If the coagulation involves a blue cycle, go to (CL);
otherwise, go to (CS).
\item[$\bullet$ (CL)] If the blue cycle coagulates with a red cycle,
and this is not the first transposition in the step
that involves a small cycle, then create a ghost
cycle; otherwise, do not create a ghost cycle.
\textit{Finish.}
\item[$\bullet$ (CS)] If a small cycle remains after the coagulation,
it is colored black.
If the coagulation involved two red cycles of size
$\ell$ and $\ell'$, create two ghost cycles of sizes $\ell$ and
$\ell'$, unless one
of these two red cycles (say of size $\ell'$) was created in the
current step, in which case create
only one ghost cycle of size $\ell$.
\textit{Finish.}\vspace*{-10pt}
\end{itemize}

\rule{\tablewidth}{0.5pt}
\legend{In addition to this description, all ghost cycles are
killed instantaneously at rate $\mu(t)$ defined
in (\ref{mu}).}
\end{table}

Let $G_i(t)$ denote the number of ghost cycles of size $i$ at
time $t$, and let $Y_i = R_i +G_i$, which counts the number of red and
ghost cycles of size $i$. Our goal is twofold. First, we
want to show that $(Y_i(t))_{i=1}^{K}$ is close in total variation
distance to $(Z_i)_{i=1}^{K}$ and second, that at time $t = t(n)$ the
probability that there is any black cycle or a ghost cycle converges to
0 as $n \to\infty$.
\begin{rmk}
Note that with our definitions,
at each step
at most one red cycle can be created,
and at most
one red cycle can disappear without being compensated by
the creation of a ghost. Furthermore these two events cannot
occur in the same step.
\end{rmk}
\begin{lemma}
\label{L:Yeq}
Assume (\ref{assumptionstep2})
as well as (\ref{assumptionstep2prime}), and let
$t = t(n)$ be as in Proposition \ref{P:step2}. Then
\[
\| (Y_i(t))_{i=1}^{K} - (Z_i)_{i=1}^{K}\|
\longrightarrow0.
\]
\end{lemma}
\begin{pf}
The idea is to observe that $Y_i$ has approximately the following dynamics:
\[
\cases{
\mbox{rate: } (x \to x+1) = \lambda, &\quad if $x\ge0$,\cr
\mbox{rate: } (x \to x-1) = i x \mu, &\quad if $x \ge1$,}
\]
and that $\lambda= \mu= k/n+o(1/n)$, so that $(Y_i)$ is approximately
a system of $M/M/\infty$ queues where the arrival rate is
$k/n$ and the departure rate of every customer is $i k/n$.
The equilibrium distribution of $(Y_i)$ is thus approximately Poisson
with parameter the ratio of the two rates, that is, $1/i$.
The number of initial customers in the queues is, by assumption
(\ref{assumptionstep2}), small enough so
that by time
$t(n)$ they are all gone, and thus the queue has reached equilibrium.

We now make this heuristics precise. To increase $Y_i$ by 1,
that is, to create a red cycle,
one needs to specify the $j$th transposition,
$1 \le j \le k-1$, of the $k$-cycle at
which it is created. The first point $x_1$ of the $k$-cycle must fall
somewhere in a blue cycle (which has probability $\theta$).
Say that $x_1 \in C_1$, with $C_1$ a blue cycle. In order to
create a cycle of size exactly $i$ at this transposition, the second point
$x_2$ must fall at either of exactly \textit{two} places within
$C_1$: either $\sigma^{i}(x_1)$ or $\sigma^{-i}(x_1)$.
However, note that if $x_2=\sigma^{-i}(x_1)$ and $|c|=k\ge3$,
then the next transposition is guaranteed to involve the newly formed cycle,
either to reabsorb it in the blue cycles,
or to turn into a black cycle through coalescence with another small
cycle or fragmentation. Either way, this newly formed cycle does not
eventually lead to an increase in $Y_i$ since by our conventions,
we do not leave a ghost in its place.
On the other hand, if $x_2=\sigma^{i}(x_1)$ then the newly formed red cycle
will stay on as a red or a ghost cycle in the next transpositions
of the application of the cycle~$c$.
Whether it stays as a ghost or a
red cycle does not change the value of~$Y_i$,
and therefore, this event leads
to a net increase of $Y_i$ by 1.
This is true for all of the first $k-2$ transpositions of the
$k$-cycle $c$, but not for the last one,
where both $x_k = \sigma^i(x_{k-1})$ and
$x_k= \sigma^{-i}(x_{k-1})$ will create a red cycle of size~$i$.
It follows from this analysis that the total rate
$\lambda(t)$ at which $Y_i$ increases by 1 satisfies
%
\begin{equation}\label{lambda+}
\lambda(t) \le\lambda^+ = \frac{k-2}{n-k+1} + \frac2{n-k+1} =
\frac{k}{n-k+1}.
\end{equation}
%
To get a lower bound, observe that for $t \le\tau$,
$\theta(t) \ge1- K(\log n)^6/n$ at the beginning of the step.
When a $k$-cycle is applied and we decompose it
into \mbox{$k-1$} elementary transpositions,
the value $\theta(t)$ for each of the transpositions
may take different successive values which we denote by
\mbox{$\theta(t,j),j=1, \ldots, k-1$}. However, note that at each such transposition,
$\theta$ can only change by at most $\pm2K/n$.
Thus it is also the case
that for all $1 \le j \le k-1$, $\theta(t,j) \ge1- 2(k-1)K(\log n)^6/n$.
Therefore,
the probability that a fragmentation of a blue cycle does not create
any small cycle is also bounded below by
%
\[
1-2(k-1)K(\log n)^6/n-2K(\log n)^6/n=1-2kK(\log n)^6/n=:\theta_-(t).
\]
It thus follows that the total rate $\lambda(t)$ is bounded below by
%
\begin{equation}\label{lambda-}
\lambda(t) \ge\theta_-^{k-1}
\biggl(\frac2n + \frac{k-2}{n}\biggr)
\ge\frac{k}n\biggl(1- 8k\frac{K(\log n)^6}n\biggr) =:
\lambda^-.
\end{equation}
Of course, by this we mean that the $Y_i(t)$ are nonnegative jump
processes whose jumps are of size $\pm1$, and that if
$\cF_t$ is the filtration generated by the entire process up to time $t$,
then
%
\begin{equation}\qquad
\lim_{h \to0^+} \frac{\PP(Y_i(t+h) = x+1 | \cF_t, Y_i(t) = x)}h =
\lambda(t) \quad\mbox{and}\quad \lambda^- \le\lambda(t) \le\lambda^+
\end{equation}
almost surely on the event $\{t \le\tau\}$.
As for negative jumps, we have that for $x \ge1$,
%
\begin{equation}\label{jump-}
\lim_{h \to0^+}\frac{\PP(Y_i(t+h) = x-1 | \cF_t, Y_i(t) = x)}h = ix
\mu(t),
\end{equation}
where $\mu(t)$ depends on the partition and satisfies the estimates
%
\begin{equation}\label{mu}
\mu^- \le\mu(t) \le\mu^+,
\end{equation}
where
%
\begin{equation}\label{mu+}
\mu^- := \frac{k}n\biggl(1- 8k\frac{K(\log n)^6}n\biggr)
\quad\mbox{and}\quad
\mu^+ = \frac{k}{n-k}.
\end{equation}
The reason for this is as follows. To decrease $Y_i$ by 1 by decreasing
$R_i$,
note that
the only way to get rid of a red cycle without creating
a ghost
is to coagulate it with a blue cycle at the $j$th transposition,
$1 \le j \le k-1$, with no other transpositions creating
small cycles.
The probability of this event is bounded above by
$ik/(n-k)$ and, with $\theta_-$ as above, bounded below by
\[
\frac{i\theta}{n}\theta_-^{k-2}+
\theta\frac{i}{n-1}\theta_-^{k-2}+
\theta\theta_-\frac{i}{n-2}\theta_-^{k-3}+
\cdots+
\theta\theta_-^{k-2}\frac{i}{n-k+1}
\geq\frac{ik}{n} \theta_-^{k-1} .
\]
Therefore, if in addition ghosts are each killed
independently with rate $\mu(t)$ as above, then (\ref{jump-}) holds.
More generally, if $1\le m\le K$ and $i_1< \cdots< i_m \le K $ are
pairwise distinct integers, then we may consider the vector
$(Y_{i_1}(t), \ldots, Y_{i_m}(t))$. If its current state is $x=(x_1,
\ldots, x_m)$, then it may make transitions to $x'=(x'_1, \ldots,
x'_m)$ where the two vectors $x$ and $x'$ differ by exactly one
coordinate (say the $j$th one) and $x_j - x'_j = \pm1$ (since only one
queue $Y_i$ can change at any time step, thanks to our coloring rules).
Also, writing $Y(t)$ for the vector $(Y_{i_1}(t), \ldots, Y_{i_m}(t))$,
we find
\[
\lim_{h \to0^+} \frac{\PP(Y(t+h) = x' | \cF_t, Y=x)}h =
\cases{
\lambda(t), &\quad if $x'_j = x_j +1$,\cr
i_j x_j \mu(t), &\quad if $x'_j = x_j -1$.}
\]
These observations show that we can compare $\{(Y_i(t\wedge\tau
)_{1\le
i \le K},t\ge0\}$ to a~system of independent Markov queues $\{
(Y^+_i(t\wedge\tau))_{1\le i \le K},t\ge0\}$ with respect to a common
filtration $\cF_t$, with no simultaneous jumps almost surely, and such
that the arrival rate of each $Y_i$ is $\lambda^+$, and the departure
rate of each client in $Y_i$ is~$i \mu^-$. We may also define a system
of queues $(Y^-_i)_{1\le i \le K}$ by accepting every new client of
$Y^+_i$ with probability $\lambda^-/\lambda^+$ and rejecting it
otherwise. Subsequently, each accepted client tries to depart at a~rate
\mbox{$\mu^+- \mu^-$}, or when it departs in $Y^+_i$, whichever comes first.
Then one can construct all three processes $(Y^-_i)_{1\le i \le K}$,
$(Y_i)_{1\le i \le K}$ and $(Y^+_i)_{1\le i \le K}$ on a~common
probability space in such a way that $Y^-_i(t) \le Y_i(t) \le Y^+_i(t)$
for all $t\le\tau$.

Note that if $(Z^+_i)_{1\le i \le K}$ denote independent Poisson random
variables with mean $\lambda^+/(i\mu^-)$, then $(Z^+_i)_{1\le i\le K}$
forms an invariant distribution for the system $(Y^+_i(t), t\ge
0)_{1\le i \le K}$. Let $(Z^+_i(t), t\ge0)_{1\le i \le K}$ denote the
system of Markov queues $Y^+_i$ started from its equilibrium
distribution $(Z^+_i)_{1\le i \le K}$. Then $(Y^+_i(t))_{1\le i \le K}$
and $(Z^+_i(t))_{1\le i \le K}$
can be coupled as usual by taking each coordinate to be equal after the
first time that they coincide. In particular, once all the initial
customers of $Y^+_i$ and of $Z^+_i(t)$ have departed (let us call $\tau
'$ this time), then the two processes $(Y^+_i)_{1\le i \le K}$ and
$(Z^+_i)_{1\le i \le K}$ are identical.

We now
check that this happens before $t=t(n)$ with high probability. It is an
easy exercise to check this for $Z^+_i(t)$ so we focus on $Y^+_i(t)$.
To see this, note that by (\ref{assumptionstep2prime}), there are no
more than $D\log(j+2)$ customers in every strip $[2^j, 2^{j+1})$
initially if $j \le\log_2 \log n$. Moreover, each customer departs
with rate at least $2^{j-1}/n$ when in this strip. Thus the time $\tau
'_j$ it takes for all initial customers of $Y^+$ in strip $[2^j,
2^{j+1})$ to depart is dominated by $(n/2^{j-1})\max_{1\le q\le D\log
(j+2)} E_q$, where $(E_q)_{q\ge1}$ is a collection of i.i.d. standard
exponential random variables. Hence
\[
\E(\tau'_j) \le\frac{n}{2^{j-2}} \bigl(\log_2 D + \log\log(j+4)\bigr).
\]
For larger strips we use the crude and obvious bound $M_j(0) \le n$ if
$j \ge
\log_2 \log n$. Moreover, each customer departs at rate $2^{j-1}/n$ with
$j\ge\break\lfloor\log_2 \log n\rfloor$. Thus, in distribution,
\[
\tau'_j \preceq\frac{n}{2^{j-2}} \max_{1\le q \le n} E_q
\]
so that $\E(\tau'_j) \le n \log n /2^{j-1}$ [we are using here
that
$\E(\max_{1\le q \le m} E_q) \le2 \log m$ for all $m$ large enough].\vadjust{\eject}
Since we obviously have $\tau' \le\sum_{j=0}^{\log_2 K+1} \tau'_j$,
we conclude
\[
\E(\tau') \le\sum_{j=0}^{\log_2
\log n} \frac{n}{2^{j-2}}\bigl(\log D + \log
\log(j+4)\bigr) + \sum_{j\ge\log_2 \log n} \frac{n \log n }{2^{j-1}}
\le
a(D) n,
\]
where $a(D)<\infty$ depends solely on $D$. By Markov's inequality and
since $t(n)/n \to\infty$, we conclude that $\tau' \le t$ with high
probability. We now claim that $(Y_i^-(t))_{1\le i \le K} =
(Y_i^+(t))_{1\le i \le K}$ with high probability. To see this, we note
that at equilibrium $\E(Z^+_i) = \lambda^+/(i \mu^-) \le2/i$. Therefore,
\begin{eqnarray*}
&&
\PP\bigl(Y^+_i(t) \neq Y^-_i(t) \mbox{ for some $1\le i\le K$}\bigr) \\
&&\qquad \le\E
\Biggl(\sum
_{i=1}^K Y^+_i(t) - Y^-_i(t);\tau' <t\Biggr) + \PP(\tau'>t) \\
&&\qquad \le\sum_{i=1}^K\frac2i\biggl\{\biggl(1- \frac{\lambda^-}{\lambda^+}\biggr) +
\biggl(1-\frac
{\mu^-}{\mu^+}\biggr)\biggr\} + \PP(\tau'>t)\\
&&\qquad \le16(k-1)\frac{K(\log n)^7}{n} + \PP(\tau'>t).
\end{eqnarray*}
Since we have already checked that $\PP(\tau'>t) \to0$ as $n \to
\infty
$, this shows that on the event $\{\tau' \le t \le\tau\}$ and $\{
Y^+_i(t) = Y^-_i(t)\mbox{ for all }1\le i \le K\}$ (an event of
probability asymptotically one),
$(Y_i(t))_{1\le i \le K}$
can be coupled to $(Z^+_i(t))_{1\le i \le K}$ which has the same law as
$(Z^+_i)_{1\le i \le K}$. Thus
%
\begin{equation}\label{YZ+}
\|(Y_i)_{i=1}^K - (Z^+_i)_{i=1}^K\| \longrightarrow0
\end{equation}
as $n \to\infty$. On the other hand, we claim that
\[
\|(Z_i)_{i=1}^K - (Z^+_i)_{i=1}^K\| \longrightarrow0
\]
also. Indeed, it is easy to see and well known that for $\alpha, \beta>0$
\[
\|\mathrm{Po}(\alpha) - \operatorname{Po}(\beta)\| \le1- \exp( - |\alpha-
\beta
|) \le|\alpha- \beta|.
\]
Since the coordinates of $Z_i$ and $Z^+_i$ are both independent Poisson
random variables but with different parameters, we find that
\begin{eqnarray*}
\|(Z_i)_{i=1}^K - (Z^+_i)_{i=1}^K\| & \le & \sum_{i=1}^K
\frac
{\lambda^+}{i \mu^-} - \frac1i\\
& \le & \sum_{i=1}^K \frac1i \biggl(\frac1{1-2(k-1)K(\log n)^6/n}
-1
\biggr)\\
& \le &\frac{4(k-1)K(\log n)^7}n \longrightarrow0
\end{eqnarray*}
as $N \to\infty$. By the
triangle
inequality and (\ref{YZ+}), this completes the proof of Lem\-ma~\ref{L:Yeq}.
\end{pf}
\begin{lemma}
\label{L:black1}
Let $t = t(n)$ be as in Proposition \ref{P:step2}. Then, with
probability tending to 1 as $n \to\infty$, $B_i(t) = 0$ for all $1
\le i \le K$.\vspace*{6pt}
\end{lemma}
%
%
\begin{pf}
Let us consider black cycles in scale $j$, that is,
those whose size~$i$ satisfies $2^j \le i <2^{j+1}$
with $j \le\log_2 K $. By assumption (\ref{assumptionstep2}),
before time $t$ the total mass of small cycles never exceeds
$2K (\log n)^6$ with high probability. Thus the rate at which a black
cycle in scale $j$ is generated by fragmentation of a red
cycle (or from another black cycle) is at most
\[
\lambda^{{B,1}}_{j} = k\frac{2K(\log n)^6}{n} \frac{2^{j+1}}n.
\]
Black cycles can also be generated directly by fragmenting a blue
cycle and subsequently fragmenting either the small cycle thus created
or some other blue cycle
in the rest of the step.
The rate at which a black fragment in scale~$j$ occurs in this fashion
is thus smaller than
\[
\lambda^{B,2}_{j} = k^2\frac{K}{n} \frac{2^{j+1}}n.
\]
Finally, one needs to deal with black cycles that arise through the
fragmentation of a blue cycle whose size at the time of the
fragmentation is between~$K$ and $2K$ (thus potentially leaving two
small cycles instead of one). Let $j' = \log_2 K$. We know that, while
$s \le\tau$, $M_{j'}(s) \le(\log n)^6/2$. In between steps, the
number of cycles in scale $j'$ cannot ever increase by more than~$2k$.
Thus the rate at which black cycles occur in this fashion at scale $j$
is at most
\[
\lambda^{(B,3)}_j=
\cases{
0, &\quad if $j < j'-1$,\vspace*{2pt}\cr
k\dfrac{K(\log n)^6}{n} \dfrac{2^{j+1}}n, &\quad if $j=j'-1$.}
\]
This combined rate is therefore smaller than $\lambda^B_j = 3 \lambda
^{B,1}_j$. Note that it may be the case that several black cycles are
produced in one step, although this number may not exceed $2k$.
On the other hand, every black cycle departs at a rate which is at least
\[
\mu^{B}_{j}= \frac{\theta}n 2^{j} \ge\frac{2^{j-1}}n
\]
since $\theta\ge1/2$ for $t \le\tau$, say. (Note that when two back
cycles coalesce, the new black cycle has an even greater departure rate
than either piece before the coalescence, so ignoring these events can
only increase stochastically the total number of black cycles.) Thus we
see that the number of black cycles in this scale is dominated by a
Markov chain $(\beta_j(s), s\ge0)$ where the rate of jumps from $x$ to
$x+2k$ is $\lambda^B_j$ and the rate of jumps from $x$ to $x-1$ is
$\mu
^B_j$, and $\beta_j(0)=0$. Speeding up time by $n/2^{j-1}$, $\beta_j$
becomes a Markov chain $\beta'_j$ whose rates are, respectively,
$\lambda'^B_j= 6k K (\log n )^6/n$ and 1, and where $\beta'_j(0)=0$. We
are interested in
\[
\PP\bigl(\beta_j(t)>0\bigr) = \PP\bigl(\beta'_j(t')>0\bigr)\qquad \mbox{where } t' = t2^{j-1}/n.
\]
Note that when there is a jump of size $2k$ (i.e., when $2k$
individuals are born) the time it takes for them to all die in this new
time-scale is a random variable $E$ which has the same distribution as
$E=\max_{1\le j \le2k} E_j$ where $(E_j)_{1\le j}$ are i.i.d. standard
exponential random variables.
Decomposing on possible birth times of individuals, and noting that $\PP
(E>x) \le2k e^{-x}$ by a~simple union bound, we see that
\begin{eqnarray*}
\PP\bigl(\beta'_j(t')>0\bigr) &=& \int_0^{t'} \lambda'^B_j \PP(E>t' -s) \,ds\\
&\le& \frac{6k K (\log n)^6 }n\int_0^\infty\PP(E>x) \,dx
\le\frac{12k^2 K (\log n)^6}n.
\end{eqnarray*}
There are $\log_2K$ possible scales to sum on, so by a union bound
the probability that there is any black cycle at time $t$ is, for large $n$,
smaller
than or equal to $ k^2 K(\log n)^8/n\to_{n\to\infty}0$.
\end{pf}

The case of ghost particles is treated as follows.
\begin{lemma}
\label{L:ghost}
Let $t = t(n)$ be as in Proposition \ref{P:step2}. Then, with
probability tending to 1 as $n \to\infty$, $G_i(t) = 0$ for all $1
\le i \le K$.\vspace*{-2pt}
\end{lemma}
\begin{pf}
Suppose a red cycle is created, and consider what happens to it the
next time it is touched. With probability at least $\theta^{k-2}$ this
will be to coagulate with a blue cycle with no other small cycle being
touched in that step, in which case this cycle is not transformed into
a ghost. However, in other cases it might become a ghost. It follows
that any given cycle in $Y_i$ is in fact a ghost with probability at most
\[
\frac{1-\theta^{k-2}}{\theta^{k-2}} \le(k-2) \frac{K (\log n)^6}n.
\]
It follows that (using the notation from Lemma \ref{L:Yeq})
\begin{eqnarray*}
\PP\bigl(G_i(t)>0 \mbox{ for some $i$}\bigr) &\le& \sum_{i=1}^K \E\bigl(G_i(t) ;
\tau
'<t\bigr) + \PP(\tau'>t)\\
&\le& \PP(\tau'>t) +\sum_{i=1}^K \frac2i\frac{(k-2)K (\log n)^6}{n}
\\
&\le& \PP(\tau'>t) + 2(k-2) \frac{K (\log n)^7}{n},
\end{eqnarray*}
which tends to 0 as $n \to\infty$.
This completes the proof of Lemma \ref{L:ghost}.\vadjust{\eject}
\end{pf}

\textit{Completion of the proof of Proposition}
\ref{P:step2}:
Since $N_i(t) = Y_i - G_i + B_i$, we get the
proposition
by combining Lemmas
\ref{L:Yeq}, \ref{L:black1} and \ref{L:ghost}.
\end{pf*}

\subsection{\texorpdfstring{Verification of (\protect\ref{assumptionstep2}) and (\protect\ref{assumptionstep2prime})}
{Verification of (8) and (9)}}
\label{subsec-verif}

In order for Proposition \ref{P:step2} to be useful, we need to show
that assumptions (\ref{assumptionstep2}) and
(\ref{assumptionstep2prime})
indeed hold with large enough probability.
This will be accomplished in Propositions
\ref{rest} and \ref{restprime} below.

Recall the variable $M_j$ [see (\ref{eq-mj})], and
let
\[
\mathcal{A}_j^s=\Bigl\{
\max_{t\in[sn \log\log n, n\log n]} M_j(t)
<n2^{-j}/(\log n)^3\Bigr\}.
\]
Recall that $K$ is the dyadic integer closest to $\lfloor n^\chi
\rfloor$.

We begin with the following lemma.
Its proof is a warm-up to the subsequent analysis.
\begin{lemma}
\label{lemcalA}
Let
\[
\mathcal{A}_\chi=\bigcap_{j=0}^{\log_2K+1} \mathcal{A}_j^{6}.
\]
Then,
\[
\PP(\mathcal{A}_\chi^\complement)\underset{\longrightarrow}{n\to
\infty}
0.
\]
\end{lemma}
\begin{pf}
It is convenient to reformulate the cycle chain as
a chain that at independent exponential times (with parameter $k$),
makes a
random transposition, where the
$\ell$th transposition is chosen uniformly at random
(if $\ell-1$ is an integer multiple of $k$), or uniformly
among those transpositions that involve the ending point of the previous
transposition and that
would result with a legitimate $k$-cycle (i.e., no repetitions
are allowed) if $\ell-1$ is not an integer multiple of $k$.

We begin with $j=0$. Note that
$M_0(0)\leq n$ and that $M_0(t)$ decreases by $1$ with rate
at least $k M_0(t)n^{-1}$
and increases, at most by $2$, with rate bounded above
by $k(1-M_0(t)/n) n^{-1}$.
In particular, by time $n\log n$, the number of increase
events is dominated by twice a Poisson variable of
parameter $k\log n$. Thus, with probability bounded
below by $1-e^{-(\log n)^2}$, at most $2(\log n)^2$ parts of size
$1$ have been born. On this event, $M_0(t)\leq2(\log n)^2+\tilde M_0(t)$
where $\tilde M_0(t)$ is a process with death only at rate
$k\tilde M_0(t)/n$. In particular,
the time of the $n-n/2(\log n)^3$th death in $\tilde M_0(t)$ is distributed
like the random variable
\[
Z_0:=\sum_{i=0}^{n-n/2(\log n)^3} \mathcal{E}_i,
\]
where the $\mathcal{E}_i$ are independent exponential
random variables of parameter
$k(n-i)/n$.
It follows that
$\E(Z_0)\sim3n \log\log n/k$ and the Chebyshev bound\vadjust{\eject}
gives, with $\zeta>0$,
\begin{eqnarray*}
\PP\bigl(Z_0>2\E(Z_0) \bigr)&\leq& \E(e^{\zeta Z_0}) e^{-2\zeta\E Z_0}
\\
&\leq& e^{-\sum_{i=0}^{n-n/2(\log n)^3}\log(1-\zeta n/k(n-i))}
e^{-6\zeta n \log\log n/k}\\
&\leq& c^{-1}e^{-c n/(\log n)^3}
\end{eqnarray*}
for an appropriate constant $c$, by choosing
$\zeta=k/2(\log n)^3$.
We thus conclude that
\[
\PP( (\mathcal{A}_0^{6/k})^\complement)\leq2e^{-(\log n)^2}.
\]

We continue on the event $\mathcal{A}_0^{6/k}$. We consider
the process
$\bar M_1(t)=M_1(t+6n\log\log n/k)$. By definition
$\bar M_1(0)\leq n/2$. The difference in the analysis of
$\bar M_1(t)$ and $M_0(t)$ lies in the fact
that now, $\bar M_1(t)$ may increase due to a~merging of two parts of
size $1$, and the departure rate is now bounded below by
$2k \bar M_1(t) n^{-1}$. Note that by time $n\log n$, the
total number of arrivals due to a merging of parts of size $1$ has mean
bounded by $n\log n \cdot k (1/(\log n)^3)^2< kn /(\log n)^6$. Repeating
the analysis concerning $M_0$, we conclude similarly that
\[
\PP( (\mathcal{A}_1^{6/k+3/k})^\complement|\mathcal{A}_0^{6/k})
\leq2e^{-(\log n)^2}.
\]

The analysis concerning $M_j(t)$ proceeds with one important
difference.
Let
$s_j=6\sum_{i=0}^{j} 2^{-i}/k $, $T_j=s_j n \log\log n$, and set
$\bar M_j(t)=M_j(t+T_{j-1})$.\vspace*{1pt} Now, $\bar M_j(t)$ can increase
due to the merging of a part of size $[2^{j-1}n,2^jn)$ with a~part
of size smaller than $2^j n$. On
$\bigcap_{i=0}^{j-1} \mathcal{A}_i^{s_i}$, this
has rate bounded above by
\[
k\frac{1}{ (\log n)^3}\cdot\frac{j}{(\log n)^3}
\leq k \frac{1}{(\log n)^5}.
\]
%
One can bound brutally the total number of such arrivals, but such a~bound is
not useful. Instead, we use the definition of the events $\mathcal{A}_i^{s_i}$,
that allow one to control the number of arrivals ``from below.''
Indeed, note that
the rate of departures $D_t$ is
bounded below by $k2^j[\bar M_j(t)-1]_+(1-1/(\log n)^2)/n$
(because the total mass below $2^j$ at
times $t\in[T_j, n\log n]$
is, on $\bigcap_{i=0}^{j-1} \mathcal{A}_i^{s_i}$, bounded above by $jn/(\log n)^3
<n/(\log n)^2$). Thus, when
$\bar M_j(t)> n2^{-j-1}/\allowbreak(\log n)^3$, the rate of departure
$D_t \gg k\frac{1}{(\log n)^5}$.
Analyzing this simple birth--death chain,
one concludes that
\[
\PP\Biggl( (\mathcal{A}_j^{s_{j}})^\complement\Big|\bigcap_{i=0}^{j-1}\mathcal
{A}_i^{s_i}\Biggr)
\leq2e^{-(\log n)^2}.
\]
Since $T_j<12 n\log\log n/k\leq6 n\log\log n$,
this completes the proof.
\end{pf}

An important corollary is the following control on the total mass
of large parts.
\begin{cor}
\label{cor-masstop}
Let $m_\chi(t)=\sum_{i>n^\chi} N_i(t)$.
Then,
\[
\lim_{n\to\infty}\PP\biggl(\min_{t\in[6n\log\log n, n\log n]}
m_\chi(t)<
n\biggl(1-\frac{1}{(\log n)^2}\biggr)\biggr)=0.
\]
\end{cor}

The next step is the following.\vspace*{5pt}
\begin{lemma}
\label{firstR}
$\!\!\!\!\!$Set $\mathcal{B}_j\,{=}\,\{\max_{t\in[k^{-1}n(\log n-\log\log n-1), n\log
n)}\!M_j(t)\,{\leq}\,(\log n)^6/2\}$. Then,
\[
\lim_{n\to\infty}
\PP\Biggl( \bigcup_{j=0}^{2 \log_2( \log n)}
\mathcal{B}_j^\complement\Biggr)=0.
\]
\end{lemma}

The proof of Lemma \ref{firstR}, while conceptually simple,
requires the introduction of some machinery
and thus is deferred to the end of this subsection.
Equipped with Lemma \ref{firstR}, we can complete the proof
of the following proposition.
\begin{prop}
\label{rest}
With notation as above,
\[
\lim_{n\to\infty}
\PP\biggl( \max_{t\in[k^{-1} n(\log n-\log\log n), n\log n]}
\max_{j=0}^{\log_2K+1} M_j(t)>(\log n)^6/2\biggr)=0.
\]
\end{prop}
\begin{pf}
Let $R=R(n)=2\log_2(\log n)$.
Because of Lemma \ref{firstR}, it is enough to consider
$M_j(t)$ for $j>R$.

We begin\vspace*{1pt} by considering $M_{R+1}(t)$.
Let ${B}_R$ denote the intersection of $\bigcap_{j=0}^R \mathcal{B}_j$
with
the complement of the
event inside the probability
in
Corollary \ref{cor-masstop}.
On the event $B_{R}$,
for $t>k^{-1}n[\log n-\log\log n-1]:= T_R$,
the rate of arrivals
due to merging of parts smaller than $2^{R}$ is bounded\vspace*{1pt}
above by
$k (2^R (\log(n))^6/n)^2$.
The rate of arrivals due to parts
larger than $2^R$ is bounded above by
$ k (2^R/n)$, and the jump is no more than 2.
Thus, the total rate
of arrival is bounded above
by $k 2^{R+1}/n$. The rate of departure on the other hand
is, due to Corollary \ref{cor-masstop}, bounded below by
$kM_{R+1}(t)2^R/n
\cdot(1-1/(\log n)^2)$. Thus, for $M_{R+1}(t)>\log n/2$, the
difference between the departure rate and the arrival rate
is bounded below by
$kM_{R+1}(t)2^R/2n$.
By definition, $M_{R+1}(T_R)\leq n2^{-R}$.
Define $T_{R+1}=T_R+ n\log n 2^{-R} $.
Let $C_{R+1}=\{\max_{t\in[T_{R+1}, n\log n]} M_{R+1}(t)<\log n\}$.
Then,
reasoning as in the proof of Lem\-ma~\ref{lemcalA}, we find that
\[
\PP(C_{R+1}^\complement|B_R)\leq e^{-(\log n)^2}.
\]
Let $B_{R+1}=B_R\cap C_{R+1}$.

One proceeds by induction. Letting
$T_{R+j}=T_{R+j-1}+ n\log n 2^{-R-j+1} $,\break
$C_{R+j}=\{\max_{t\in[T_{R+j}, n\log n]} M_{R+j}(t)<\log n\}$ and
$B_{R+j} = B_{R+j-1} \cap C_{R+j}$, we
obtain from the same analysis that
for $j=1,\ldots,\log(K)+1$,
\[
\PP(C_{R+j+1}^\complement|B_{R+j})\leq e^{-(\log n)^2}.
\]
Thus,
$\PP(B_{R+\log(K)+1}^\complement)\leq\PP(B_R^\complement) +(\log n)e^{-(\log n)^2} \to_{n\to\infty} 0$,
while\break
$T_{R+\log(K)+1}\leq
k^{-1}n[\log n- \log\log n -1+2^{-R} \log n \sum_{j\ge1} 2^{-j}]$.
This
completes the proof, since $2^{R}=(\log n)^{2}$.
\end{pf}

\subsection{\texorpdfstring{Proof of Lemma \protect\ref{firstR}}{Proof of Lemma 10}}

While a proof could be given in the spirit of the proof
of Lemma \ref{lemcalA}, we prefer to present a conceptually simple
proof based on comparison with the random $k$-regular hypergraph. This
coupling is analog to the usual coupling with an Erd\H{o}s--R\'{e}nyi
random graph (see, e.g., \cite{BD} and~\cite{schramm}).
Toward this end, we need the following definitions.
\begin{defi}
\label{def-hyper}
A \textit{$k$-regular hypergraph} is a pair $G=(V,H)$ where $V$ is
a (finite) collection of vertices, and $H$ is a collection
of subsets of $V$ of size~$k$. The
\textit{random} hypergraph $G_k(n,p)$ is defined as the
hypergraph consisting of $V=\{1,\ldots,n\}$,
with each subset $h$ of $V$ with $|h|=k$ taken independently
to belong to $G_k(n,p)$ with probability $p$.
\end{defi}

Let $G_t$ denote the random $k$-hypergraph obtained by
taking $V=\{1,\ldots,n\}$ and taking $H$ to consist of the
$k$-hyperedges corresponding
to the $k$-cycles $\gamma_1,\ldots,\gamma_{N_t}$
of the random walk $\pi_t$.
It is immediate to check that $G_t$ is distributed like
$G_k(n,p_t)$ with
\[
p_t=1-\exp\biggl(-\frac{t}{{n\choose k}} \biggr)
\sim\frac{k! t}{n^k}.
\]
\begin{defi}
\label{def-hyper1}
A \textit{$k$-hypertree} with $h$ hyperedges in a $k$-regular
hypergraph $G$ is a connected component of $G$ with
$i=(k-1)h+1$ vertices.
\end{defi}

(Pictorially, a $k$-hypertree corresponds to a standard
tree with hyperedges, where any two hyperedges have at most
one vertex in common.)
$k$-hypertrees can be easily enumerated, as in the following, which is
Lemma 1 of \cite{KL}.
\begin{lemma} The number of $k$-hypertrees with $i$ (labeled) vertices
is
%
\begin{equation}
\label{eq-hypertreenum}
\frac{[(k-1)h]!i^{h-1}}{h! ((k-1)!)^h},\qquad h\geq0,
\end{equation}
where $h$ is the number of hyperedges and thus
$i=(k-1)h+1$.
\end{lemma}

The next lemma controls the number of $k$-hypertrees with
a prescribed number of edges in $G_t$.\vspace*{2pt}
\begin{lemma}\label{cont-tree}
Let
\begin{eqnarray*}
\mathcal{D}_{t,h}&=&\{\mbox{\# of $k$-hypertrees with \mbox{$\le$}$ h$ hyperedges in
$G_t$}\\
&&\hspace*{76.5pt}\mbox{is not larger than $(\log n)^{1.1}$}\}.
\end{eqnarray*}
Then,
%
\begin{equation}
\label{eq-boundnew}
\PP\biggl(\bigcap_{t>(n/k)[\log n-\log\log n-1
]}
\mathcal{D}_{t,(\log n)^2}\biggr)
\underset{\longrightarrow}{n\to\infty} 1 .
\end{equation}\vspace*{2pt}
\end{lemma}
\begin{pf}
Let $t_0=k^{-1}n[\log n-\log\log n-1]$ and $h_0=(\log n)^2$.
By monotonicity, it is enough to check that
%
\begin{equation}
\label{eq-boundnew1}
\PP( \mathcal{D}_{t_0,h_0})
\underset{\longrightarrow}{n\to\infty} 1 .
\end{equation}
Note that, with $i=(k-1)h+1$,
and adopting as a convention $h \log h = 0$ when $h=0$,
%
\begin{eqnarray}\label{ET_h}\quad
\PP(
\mathcal{D}_{t_0,h_0}^\complement)
&\leq&
\sum_{h=0}^{(\log n)^2}
\frac{
\E(\mbox{$\#$ of $k$-hypertrees
with $ h$ hyperedges in $G_{t_0}$)}}{(\log n)^{1.1}} \nonumber\\
&\leq&
\frac{1}{(\log n)^{1.1}}
\sum_{h=0}^{(\log n)^2}
\pmatrix{n \cr i}
\frac{( (k-1)h)! i^{h-1}}
{h! ( (k-1)!)^h}
p_{t_0}^h
(1-p_{t_0})^{{i \choose k}-h+ i{n-i \choose k-1} }\\
&\leq&
C_k
\sum_{h=0}^{(\log n)^2}
(\log n)^{i+h-1.1} e^{-(k-1)h(\log n-\log h(k-1))}
\underset{\longrightarrow}{n\to\infty}0 .\nonumber
\end{eqnarray}
[Indeed recall that if $T$ is a subset of $\{1, \ldots, n\}$ comprising
$i$ elements, then disconnecting $T$ from the rest of $\{1, \ldots, n\}
$ requires closing exactly
${i \choose1}{n-i \choose k-1} + {i \choose2} {n-i \choose k-2} +
\cdots+ {i \choose k-1} {n-i \choose1} \ge i {n-i \choose k-1}$ hyperedges,
while
${i \choose k} - h$
is the number of hyperedges that need to be closed inside $T$ for it to
be a~hypertree.]
\end{pf}

We can now provide the following proof:\vspace*{2pt}
\begin{pf*}{Proof of Lemma \ref{firstR}}
At time $t$, $N_i(t)$ consists of cycles that have been obtained
from the coagulation of cycles that have never fragmented
during the evolution by time $t$, denoted $N_i^c(t)$,
and of cycles that have been
obtained from cycles that have fragmented and created a part
of size less than or equal to $i$, denoted $N_i^f(t)$.
Note that $N_i^c(t)$ is dominated above by the number of
$k$-hypertrees with $h$ edges in $G_t$, where $i=(k-1)h +1$.
By Lemma \ref{cont-tree}, this is bounded
above by $(\log n)^{1.1}$ with high probability
for all $i\le(\log n)^2$.
On the other hand,
the rate of creation by fragmentation of cycles of size $i$
is bounded above by $4k/n$, and hence by time $n\log n$,
with probability approaching $1$
no more than $(\log n)^{1.1}$ cycles of size $i$ have been
created, for all $i\leq(\log n)^2$.
We thus conclude that with probability tending to $1$,
we have, with
$t_0=k^{-1}n[\log n-\log\log n-1]$,
%
\[
\max_{i\leq(\log n)^{2}} \max_{
t\in[t_0,n\log n]}
N_i^f(t)
\leq(\log n)^{3.1}.
\]
This yields the lemma, since for
$j\leq2 \log_2(\log n)$,
%
\[
M_j(t)\leq(\log n)^{2} \max_{i\leq(\log n)^{2}} N_i(t).
\]
\upqed\end{pf*}

\subsection{\texorpdfstring{Proof of (\protect\ref{assumptionstep2prime})}{Proof of (9)}}

We now prove that at time $\tm= (1/k) n \log n$,
the assumption (\ref{assumptionstep2prime}) [with $M_j(0)$ replaced
by $M_j(\tm)$] is satisfied, with
high probability.
\begin{prop} \label{restprime} For every $\eps>0$ there
exist $D=D(\eps)>0$ and $n_0=n_0(\eps)$
such that for $n>n_0$,
\[
\PP\bigl(M_j(\tm) \le D \log(2+j), j=0,1,\ldots,\log_2
\log n + 1\bigr)\geq(1-\eps).
\]
\end{prop}
\begin{pf}
Consider first the time $u= \frac1k (n\log n - n \log\log n)$.
\begin{lemma}\label{restprime-u}
With probability approaching 1 as $n\to\infty$, we have $M_j(u)\le
2^{j+4} \log n $ for all $0\le j \le\log_2 n$.
\end{lemma}
\begin{pf}
As in the proof of Lemma \ref{firstR},
split $M_j(t)$ into two components $M_j^f(t)$ and $M_j^c(t)$.
Note that the rate at which a fragment of size
less than $2^{j+1}$ is produced is smaller
than $2^{j+2}k/n$, so for any $w\le(1/k) n
\log n$, $M^f_j(w) \le\Poi(2^{j+2} \log n)$.
The probability that such a Poisson random variable
is more than twice its expectation
is (by standard
large deviation bounds)
smaller than $n^{-\alpha}$ for some $\alpha>0$, so summing over $\log
_2 \log n$ values of
$j$ we easily obtain that with high probability, $M_j^f(u) \le2^{j+3}
\log
n$ for all $0\le j \le\log_2 \log n$.

It remains to show that $M^c_j(u) \le\log n$ for all
$0\le j \le\log_2 \log n$ with
high probability. To deal with this part, note that if $T_h$ denotes
the number of hypertrees with $h$ hyperedges in $G_u$, then
$N_i^c(u) \le T_h$ where $i=1+h(k-1)$ is the number of vertices.
Reasoning as in (\ref{ET_h}), we compute after simplifications
[recalling that $u = (1/k)(n \log n - n \log\log n)$
and $i =1+h(k-1)$], for $h\ge0$
%
\begin{eqnarray} \label{boundET_h}
\E(T_h) & = &\pmatrix{n \cr i} \frac{(i-1)! i^{h-1}}{h! ((k-1)!)^h} p_u^h
(1-p_u)^{{i \choose k} - h + i{n-i \choose k-1}} \nonumber\\[-8pt]\\[-8pt]
&\le&\frac{n (\log n)^h }{h! i }(1-p_u)^{i{n-i \choose k-1}}
\le\frac{n^{1-i}(\log n)^{1+hk}}{h!i}.\nonumber
\end{eqnarray}
Thus summing over $i=2$ to $i = \lceil\log n \rceil$,
we conclude by Markov's inequality that $M^c_j(u) =0$ for all
$1\le j \le\log_2
\log n$ with high probability.
For $i=1$ or $h=0$, we get from (\ref{boundET_h})
\[
\E(T_0) \le\log n.
\]
Computing the variance is easy: writing $T_0 = \sum_{v \in V} \mathbf
{1}_{\{v\ \mathrm{ is}\ \mathrm{isolated}\}}$, we get
\[
\var(T_0) \le\E(T_0) + \sum_{v \neq w} \cov\bigl(\mathbf{1}_{\{v\
\mathrm{is}\ \mathrm{isolated}\}}, \mathbf{1}_{\{w\ \mathrm{is}\ \mathrm{isolated}\}}\bigr).
\]
But note that
\[
\PP(v \mbox{ is isolated}, w \mbox{ is isolated}) = \frac{\PP(v
\mbox{ is
isolated})^2}{1-p_u},
\]
so
\[
\var(T_0) \le\E(T_0) + \E(T_0)^2\biggl(\frac1{1-p_u} -1\biggr)\le
\E
(T_h) + o(1).
\]
Thus by Chebyshev's inequality, $\PP(M_0^c(u)> 2 \log
n) \to0$ as $n \to\infty$. This proves the lemma.\vspace*{-2pt}
\end{pf}

With this lemma we now complete the proof of Proposition \ref{restprime}.
We compare $(M_j(t),t\ge u)$ to independent queues as follows. By
Proposition \ref{rest}, on an event of high probability, during the
interval $[u,\tm]$ the rate at which some two cycles of size
smaller than $\log n$ coagulate is smaller than $O(((\log n)^7/n)^2)$,
so the probability that this happens during this interval of time is $o(1)$.
Likewise, the rate at which some
cluster smaller than $\log n$ will fragment is at most $k (\log
n)^{14} /n^2$,
so the probability that this happens during the interval $[u, \tm]$ is
$o(1)$. Now, aside from rejecting any $k$-cycle that would create such
a transition, the only possible transition for $M_j$ are increases by 1
(through the fragmentation of a component larger than $2\log n$) and
decreases by 1 (through coagulation with cycle larger than $\log n$).
The respective rates of these transitions is, as in (\ref{lambda+}), at
most $2^j \lambda^+ = 2^jk/(n-k)$, and at least $\nu=
2^j(k/n)(1-(\log
n)^3/n))$ as in (\ref{mu+}).
This can be compared to a queue where both the departure rate and the
arrival rate are equal to~$\lambda^+$, say $\bar M_j(t)$. The
difference between $M_j(t)$ and $\bar M_j(t)$ is that some of the
customers having left in $\bar M_j(t)$ might not have left yet in
$M_j(t)$. Excluding the initial customers, a
total of $\Poi(2^j \log\log
n )$ customers arrive in the queue $\bar M_j(t)$ during the interval
$[u, \tm]$, so the probability that any one of those customers has not
yet left by time $\tm$ in $M_j(t)$ given that it did leave in $\bar
M_j(t)$ is no more than $\lambda^+/\nu-1 = O((\log
n)^3/n)$, where the constants implicit in $O(\cdot)$ do not depend on
$j$ or $n$. Thus with probability greater than $1- O(2^j \log\log n
(\log n)^3 /
n)$, there is no difference between~$M_j(\tm)$ and $\bar M_j(\tm)$. Moreover,
%
\begin{equation}\label{boundMj1}
\bar M_j(\tm) \preceq\Poi(1) + R_j,
\end{equation}
where $R_j$ is the total number of initial customers customers that
have not departed yet by time $\tm$. Using Lemma \ref{restprime-u},
%
\begin{equation}\label{boundRj}
\{R_j>0\} \subset\biggl\{\frac1{\lambda^+}\max_{1\le q \le2^{j+4} \log
_2 n
} E_q
<\tm\biggr\},
\end{equation}
where $(E_q, q\ge1)$ is a collection of i.i.d. standard exponential
random variables. Using the independence of the queues $\bar M_j(t)$,
in combination with (\ref{boundMj1}) and (\ref{boundRj}) as well as standard
large deviations for Poisson random variables, the proposition follows
immediately.
\end{pf}

\subsection{Conclusion: Small cycles}
\label{subsec-conclusionrev}
Combining Propositions \ref{P:step2} and \ref{rest},
and using the notation introduced in the beginning of this section,
we have proved the following. Fix $\eps>0$. Then there is a $c_{\eps
,k}>0$ such that
with $t=t(n)=k^{-1}n\log n+c_{\eps,k} n$, and
all large $n$,
%
\begin{equation}\label{conclsmall0}
\| (N_i(t))_{i=1}^{K} - (Z_i)_{i=1}^{K}\|<
\eps.
\end{equation}

We now deduce the following:\vspace*{2pt}
\begin{prop}
\label{prop-small}
Fix $\eps>0$. Then there is a $c_{\eps,k}>0$ such that
with $t=t(n)=k^{-1}n\log n+c_{\eps,k} n$, and
all large $n$,
%
\begin{equation}\label{conclsmall1}
\| (N_i(t))_{i=1}^{K} - (N_i)_{i=1}^{K}\|<
\eps,
\end{equation}
where $(N_i)_{1\le i \le n}$ is the cycle distribution of a random
permutation sampled according to the invariant distribution
$\mu$.\vspace*{2pt}
\end{prop}
\begin{pf}
By (\ref{conclsmall0}) and the triangle
inequality,
all that is needed is to show that
%
\begin{equation}\label{conclsmall2}
\| (Z_i)_{i=1}^{K} - (N_i)_{i=1}^{K}\| \to0.
\end{equation}
Whenever $k$ is even, and thus
$\mu$ is uniform on $\cS_n$, (\ref{conclsmall2})
is a classical result of Diaconis--Pitman and of Barbour, with
explicit upper bound of $4K/n$ (see~\cite{barbour} or
the discussions around \cite{ArratiaTavare}, Theorem 2,
and \cite{abt}, Theorem 4.18).

In case $k$ is odd, $\mu$ is uniform on $\cA_n$. A sample $\gamma$
from $\mu$ can be
obtained from a sample $\gamma'$
of the uniform measure on $\cS_n$ using the following
procedure. If $\gamma'$ is even, take $\gamma=\gamma'$, otherwise
let $\gamma=\pi\circ\gamma'$ where $\pi$ is some fixed transposition
[say $(12)$].
The probability that the collection of small cycles in $\gamma$ differs
from the corresponding one in $\gamma'$ is bounded above
by $4 K/n\to0$, which completes the proof.
\end{pf}

\section{Large cycles and Schramm's coupling}
\label{sec-schramm}

Fix $\eps>0$ and $\chi\in(7/8,1)$.
Recall that $K$ is the closest dyadic integer
to $\lfloor n^\chi\rfloor$ and that a cycle is called
small if its size is smaller than $K$.
For $n$ large, let $t= t(n)= k^{-1}n \log n+c_{\eps,k}n$.
We know by the previous section (see Proposition
\ref{prop-small})
that at this time, for $n$ large,
the distribution of the small cycles of
the permutation $\pi_t$ is
arbitrarily close (variational distance smaller
than $\eps$) to that of a (uniformly chosen)
random permutation $\pi'$. Therefore we can find a coupling of
$\pi:=\pi_t$ and $\pi'$ in such a way that
%
\begin{equation} \label{ic}
\mathbb{P}(\mbox{the small cycles of $\pi$ and $\pi'$ are
identical}) \ge1-\eps.
\end{equation}
We can now provide the following proof:
\begin{pf*}{Proof of Theorem \ref{T:mix}}
We will construct an evolution of $\pi'$, denoted~$\pi'_s$,
that follows the random
$k$-cycle dynamic (and hence, $\pi'_s$ has cycle structure whose law
coincides
with the law of the
cycle structure of a uniformly chosen permutation, at all times).
The idea is that with
small cycles being the hardest to mix,
coupling $\pi_{t+s}$ and $\pi'_s$ will now take very little time.
To prove this, we describe a modified version of the Schramm coupling
introduced in \cite{schramm},
which has the additional
property that it is difficult to create small unmatched pieces.

To describe this coupling,
we will need some notation from \cite{schramm}. Let $\Omega_n$ be the
set of discrete partitions of unity
\[
\Omega_n\,{=}\,\Biggl\{\!(x_1\,{\ge}\,\cdots\,{\ge}\,x_n)\dvtx x_i\,{\in}\,\{0/n, \ldots, n/n\}
\mbox{ for all } 1\,{\le}\,i\,{\le}\,n\mbox{,  and } \sum_{i=1}^n x_i\,{=}\,1\!\Biggr\}.
\]
We identify the cycle count of $\pi_t$ with a vector $Y_t \in\Omega
_n$. We thus want to describe a coupling between two processes $Y_t$
and $Z_t$ taking their values in $\Omega_n$ and started from some
arbitrary initial states. The coupling will be described by a joint
Markovian evolution of $(Y_t, Z_t)$.

We now begin by describing the construction
of a random transposition. For $x\in(0,1)$, let $\{x\}_n$ denote the
smallest element of $\{1/n,\ldots,n/n\}$ not smaller than~$x$.
Let $\tilde u, \tilde v$ be two random points uniformly
distributed in $(0,1)$,
set $u=\{\tilde u\}_n, v=\{\tilde v\}_n$
and condition them so that $u\neq v$.
Note that $u,v$ are
both uniformly distributed on $\{1/n,\ldots,n/n\}$.
If we focus for one moment on the marginal evolution of $(Y_t)$, then
applying one transposition to $Y_t$ can be realized by
associating to $Y_t \in\Omega_n$ a tiling of the semi-open
interval $(0,1]$ where each tile is equally semi-open and there
is exactly one tile for each nonzero coordinate of $Y_t$. (The order
in which those tiles are put down may be chosen arbitrarily and does
not matter for the moment.)
If $u$ and $v$ fall in different tiles then we merge the two tiles
together and get a new element of~$\Omega_n$ by sorting in decreasing
order the size of the tiles. If $u$ and $v$ fall in the same tile then
we use the location of $v$ to split that tile into two parts: one that
is to the left of $v$, and one that is to its right (we keep the same
semi-open convention for every tile). This procedure works because,
conditionally on falling in the same tile $C$ as $u$, then $v$ is
equally likely to be on any point of $C \cap\{1/n, \ldots, n/n\}$
distinct from $v$, which is the same fragmenting rule as explained at
the beginning of the proof of Proposition~\ref{P:step2}.

We now explain how to construct one step of the joint evolution.
If $Y,Z \in\Omega_n$ are two unit discrete partitions, then we can
differentiate between the entries that are matched and those that are
unmatched; two entries from $Y$ and $Z$ are matched if they are of
identical size. Our goal will be to create as many matched parts as possible.
Let $Q$ be the total mass of the unmatched parts.
When putting down the tilings associated with $Y$ and $Z$ we will do so in
such a way that all matched parts are at the
right of the interval
$(0,1]$ and the unmatched parts occupy the left part of the interval,
as in Figure~\ref{Fig:coupl1}.
If $u$ falls into the matched parts,
we do not change the coupling beyond that described
in \cite{schramm};
that is, if $v$ falls in the same component as $u$ we make the same
fragmentation in both copies, while otherwise we make the corresponding
coalescence.
The difference occurs if $u$ falls\vspace*{2pt} in the unmatched parts.
Let $y$ and $z$ be the respective components of $Y$ and $Z$ where $u$ falls,
and let $\hat Y, \hat Z$ be the reordering of $Y,Z$
in which these components
have been put to the left of the interval $(0,1]$.
%
\begin{figure}

\includegraphics{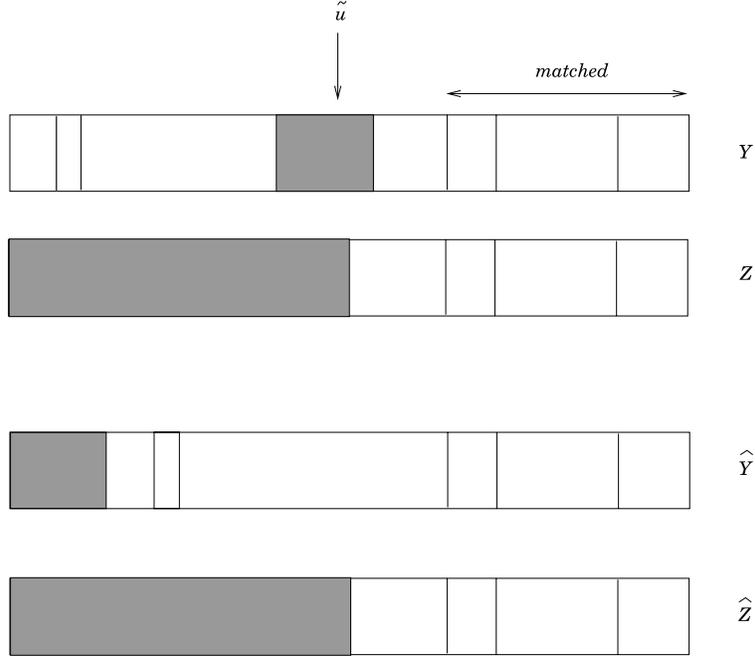}

\caption{First step of the coupling. A point $\tilde
u$ is uniformly chosen on $(0,1)$ and picks a~part
in $Y$ and~$Z$, which are then rearranged into $\hat Y, \hat Z$.}
\label{Fig:coupl1}
\end{figure}
Let $a=|y|$ and let $b=|z|$ be the respective
lengths of the pieces selected with $u$,
and assume without loss of generality that $a<b$.
Further rearrange, if needed, $y$ and $z$ so that after the rearrangement,
\mbox{$|u|=1/n$}. Because $v\neq u$, necessarily $v>1/n$ (and is uniformly
distributed on the set $\{2/n,\ldots,n/n\}$).
The point $v$ designates a~size-biased sample from the partition
$\hat Y$ and we will construct another point~$v'$, which will also be
uniformly\vspace*{1pt} distributed on $\{2/n,\ldots,n/n\}$,
to similarly select a size-biased sample
from $\hat Z$. However, while in the coupling of \cite{schramm} one takes
$v=v'$, here we do not take them equal and apply to $v$ a~measure-preserving
map $\Phi$, defined as follows.
Define the function
%
\begin{equation}\label{Phi}
\Phi(x)=
\cases{
x, &\quad if $x>b$ or if $1/n\leq x\le\gamma_n +1/n$, \cr
x - \gamma_n, &\quad if $a<x\le b$,\cr
x + b-a, &\quad if $\gamma_n +1/n<x\le a$,}
\end{equation}
%
where $\gamma_n:=\{(a-1/n)/2\}_n$. See Figure \ref{Fig:coupl2} for
description of $\Phi$.
Note that $\Phi$ is a measure-preserving map and hence $\tilde v':=
\Phi
(\tilde
v)$
is uniformly distributed on $(0,1)$. Define $v'=\{\tilde v'\}_n$.
With $u,v$ and $v'$ selected,
the rest of the algorithm is unchanged, that is, we make the
corresponding coagulations and fragmentations.

\begin{figure}

\includegraphics{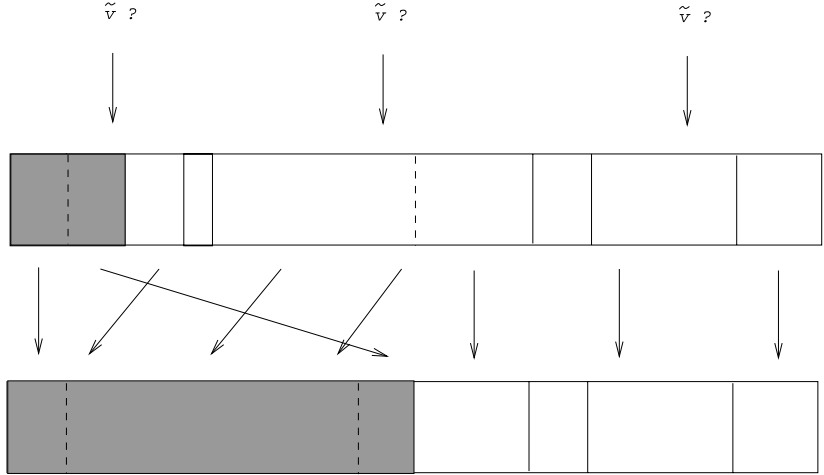}

\caption{A second point $\tilde
v$ is chosen uniformly in (0,1) and serves as a
second size-biased pick for~$\hat Y$. $\tilde
v$ is mapped to $\tilde
v'=\Phi(\tilde
v)$ which gives a second size-biased pick for $\hat Z$.}
\label{Fig:coupl2}
\end{figure}

This coupling has a number of remarkable properties which we summarize below.
Essentially,
the total number of unmatched entries can only decrease, and
furthermore it is very difficult to create small unmatched entries, as
the smallest unmatched entry can only become smaller by a factor of at
most~2.\looseness=1

In what follows, we often speak of the ``unmatched entries'' between
two permutations, meaning that we associate to these permutations
elements of~$\Omega_n$ and identify matched parts in $\Omega_n$
with matched cycles in the permutations. The translation between the two
involves a factor $n$ concerning the size of the parts, and in all places
it should be clear from the context whether we discuss parts in $\Omega
_n$ or
cycles of partitions.
\begin{lemma}\label{L:coupling}
Let $U$ be the size of the
smallest unmatched entry in two partitions $Y,Z \in\Omega_n$,
let $Y',Z'$ be the corresponding partitions after one transposition of
the coupling and let $U'$ be the size of the smallest unmatched entry
in $Y',Z'$. Assume that $2^j\le U < 2^{j+1}$ for some $j \ge0$.
Then it is always the case that $U'\ge U-\{U/2\}_n$, and moreover,
\[
\PP( U' \le2^j) \le2^{j+2}/n.
\]
Finally, the number of unmatched parts may only decrease.
\end{lemma}
\begin{rmk}
$\!\!\!$Since $U'\,{\ge}\,U\,{-}\,\{U/2\}_n$, it holds in particular that $U'\,{\ge}\,2^{j-1}$.
\end{rmk}
\begin{pf*}{Proof of Lemma \ref{L:coupling}}
That the number of unmatched entries can only decrease is similar to
the proof of Lemma 3.1 in \cite{schramm}. (In fact it is simpler here,
since that lemma requires looking at the total number of unmatched
entries of size greater than $\eps$. Since in our discrete setup no
entry can be smaller than $\eps= 1/n$ we do not have to take this
precaution.)
We continue to denote by
$M_j$ the total number of parts in the range $[2^j, 2^{j+1})/n$.
The only case that $U$ can decrease is if there is a fragmentation of an
unmatched entry, since matched entries must fragment in exactly the
same way.
Now, note that the coupling is such that when an
unmatched entry is selected and is fragmented, then all subsequent pieces
are either greater or equal to $a-\{a/2\}_n$
(where $a$ is the size of the smaller
of the two selected
unmatched entries), or are matched. Moreover, for such a fragmentation
to occur, one must select the lowest
unmatched entry (this has probability at most $M_j2^{j+1}/n$,
since there may be several
unmatched entries with size~$U$), and then fragment it, which has
probability at most $2^{j+1}/n$, and thus $\PP(U'<U) \le4 M_j 4^j/n^2$.
Since $M_j 2^j\leq n$, this completes the proof.
\end{pf*}

We have described the basic step of a (random) transposition in the coupling.
The step corresponding to a random
$k$-cycle $\gamma= (\gamma_1,\gamma_2,\ldots,\gamma_k)$
is obtained by taking $u_1=\gamma_1$, generating $v,v'$ as in the coupling
above (corresponding to the choice of $\gamma_2$), rearranging and taking
$u_2$ to correspond to the location of $v,v'$ after the rearrangement, drawing
new $v,v'$ (corresponding to $\gamma_3$) and so on.
In doing so, we are disregarding the constraint that no repetitions are present
in $\gamma$. However, as it turns out, we will
be interested in an evolution lasting
at most
%
\begin{equation}
\label{eq-070110}
\Delta:= n^{5/8} \log n ,
\end{equation}
%
and the expected number of times
that a violation of this constraint occurs
during this time is bounded by
$ 2\Delta k^2/n$, which
converges to $0$ as $n\to\infty$. Hence, we can in what follows disregard
this violation of the constraint.

Now, start with two configurations $Y_0,Z_0$ such that $Z_0$
is the element of $\Omega_n$ associated with a random uniform permutation.
Assume also that initially, the small parts of $Y_0$ and $Z_0$
(i.e., those that are smaller than $K$, the closest dyadic integer to
$\lfloor n^\chi\rfloor$),
are exactly identical, and that they have the same parity.
As we will now see, at time $\Delta$,
$\pi_{t+\Delta}$ and $\pi'_\Delta$ will be coupled,
with high probability.
Note also that, since initially all the parts that are smaller than $K$
are matched, the initial number of unmatched entries cannot exceed $n/K
\le n^{1/8}$, and this may only decrease with time by Lemma~\ref{L:coupling}.
\begin{lemma} \label{L:mass small}
In the next $\Delta$ units of time, the random permutation $\pi'_s$
never has more than
a fraction $ n^{-1/8} (\log n)^6$ of the total mass in parts
smaller than $n^{7/8}$, with high probability.
\end{lemma}
\begin{pf}
The proof is the same as that of Proposition \ref{rest}, only simpler because
the initial number of small clusters is within the required range.
We omit further details. [This can also be seen by computing the
probability that a given
uniform permutation $\pi'$ has more than a fraction
$n^{-1/8} (\log n)^6$
of the total mass in parts smaller than $n^{7/8}$,
and summing over $\Poi(\Delta)$ steps.]
\end{pf}
%
%
\begin{lemma} \label{L:unmatched large}
In the next $\Delta$ units of time, every unmatched part
of the permutations
is greater than
or equal to $n^{3/4}/2$, with high probability.
\end{lemma}
\begin{pf}
Recall that the total number of unmatched parts
can never increase. Suppose\vspace*{1pt} the smallest unmatched part
at time $s$
is of scale $j$ (i.e., of size in $[2^j, 2^{j+1})$),
and let $j=U(s)$ be this scale. Then, when touching this
part,
the smallest scale it could go to is $j-1$,
by the properties of the coupling (see Lem\-ma~\ref{L:coupling}).
This happens with probability at most $2^{j+2}/n$.
On the other hand, with the complementary probability, this part
experiences a coagulation. And with reasonable probability, what it
coagulates with is larger than itself, so that it will jump to scale
$j+1$
or larger. To compute this probability, note that since this
is the smallest unmatched part,
all smaller parts
are matched and thus have a total mass controlled by
Lemma~\ref{L:mass small}.
In particular, on an event of
high probability, this fraction of the total mass is at most
$q:= n^{-1/8} (\log n)^6$. It follows that with probability
at least $1-q$, the part
jumps to scale at least
$j+1$, and with probability at most
$r_j:=2^{j+1} /n$,
to scale $j-1$. Now, when this part
jumps to scale at least
$j+1$, this does not necessarily mean that the
\textit{smallest} unmatched part
is in scale at least
$j+1$, since there may be several small unmatched parts
in scale $j$. However, there can never be more than $2n^{1/8}$ such parts.
If an unmatched piece in scale $j$ is touched, we declare it
a success if it moves to scale $j+1$
(which has probability at least $1-q$, given that it is touched)
and a failure if
it goes to scale $j-1$
(which has probability at most $r_j$).
If $2n^{1/8}$ successes occur before any failure occurs at scale~$j$,
we say that a \textit{good success} has occurred, and then
we know that no unmatched cycle can exist at scale smaller than
$j$.
Call the complement of a good success a \textit{potential failure} (which
thus includes the cases of both a real failure and a success which is
not good).
The probability of a potential failure at scale $j$
is at most $2n^{1/8}r_j/(1-q+r_j)$, which is
bounded above by $p_j=6n^{1/8}2^j/n$.

Let $\{s_i\}_{i\geq0}$ be the times at which the smallest unmatched
part changes scale, with $s_0$ being the first time the smallest
unmatched part
is of scale $j_0$ where \mbox{$2^{j_0}=n^{5/6}$}.
Let $\{U_i\}$ denote the scale of the smallest unmatched
part at time $s_i$, and let $j_1$ be such that
$2^{j_1}=n^{3/4}/2$. Introduce a birth--death chain on
the integers, denoted $v_n$, such that $v_0=j_0$ and
%
\begin{equation}
\label{may27}
\PP(v_{n+1}=j-1|v_n=j)=
\cases{
1, &\quad if $j=j_0$,\cr
0, &\quad if $j=j_1$,\cr
p_j, &\quad otherwise,}
\end{equation}
and
%
\begin{equation}\quad
\label{may27a}
\PP(v_{n+1}=j+1|v_n=j)=
\cases{
1-\PP(v_{n+1}=j-1|v_n=j), &\quad$j>j_1$,\cr
0, &\quad$j=j_1$.}
\end{equation}
%
Set
$\tau_j=\min\{n>0\dvtx v_n=j\}$,
and
an analysis of the birth--death chain defined by
(\ref{may27}) and (\ref{may27a}) gives that
\[
\PP^{j_0}(\tau_{j_1}<\tau_{j_0}) = \frac1{\sum_{j=j_1+1}^{j_0}
\prod
_{m=j}^{j_0-1} (({1-p_m})/{p_m})} \leq\prod_{j=j_1+1}^{j_0-1}
\frac{p_j}{1-p_j}
\]
(see, e.g., Theorem (3.7) in Chapter 5 of \cite{durrett}).
Thus
$ \PP^{j_0}(\tau_{j_1}<\tau_{j_0})$ decays as an exponential in
$(\log n)^2$.
Therefore, since
$\PP(v_{2k\Delta}=j_1)\leq2k\Delta\PP^{j_0}(\tau_{j_1}<\tau_{j_0})$,
it follows that
$\PP(v_{2k\Delta}=j_1)\to0$ as $n\to\infty$.
On the other hand, between times~$t$ and $t+\Delta$, the process $\{
U_i\}_{i\ge1}$ may have made at most $2k\Delta$ moves
with overwhelming probability. This implies that $U_i \ge j_1$ with
high probability throughout $[t, t+\Delta]$.
\end{pf}

\textit{End of the proof of Theorem} \ref{T:mix}.
We now are going to prove that, after
$\Delta= n^{5/8}\log n$ steps, there are no more unmatched parts
with high probability. The basic idea is that, on the one hand,
the number of unmatched parts
may never increase, and on the other hand, it does decrease
frequently enough. Since each unmatched part
is greater than $n^{3/4}/2$ during this time, any given pair of unmatched
parts
is merging at rate roughly $n^{-1/2}$.
There are initially no more than $2n^{1/8}$ unmatched parts,
so after $n^{5/8} \log n=\Delta$ steps, no more unmatched
part
remains with high probability.

To be precise, assume that there are $L$ unmatched parts.
Let $T_L$ be the time to decrease the number of unmatched parts
from $L$ to $L-2$. Observe that, for parity reasons ($\pi$ and $\pi'$
must have the same parity of number of parts
at all times), $L$ is always even. Note also that $L=2$ is impossible,
so $L$ is at least 4. Assume to start with that both copies have at
least 2 unmatched
parts. Then, at rate greater than $n^{-1/4}/2$ we pick an unmatched
part in the first point $u_1$ for the $k$-cycle.
Since there are at least 2 unmatched parts
in each copy, let $R$ be the interval of $(0,1)$ corresponding to a
second unmatched
part in the copy that contains the larger of the two selected ones.
Then $|R| > n^{-1/4}/2$, and moreover when $v$ falls in $R$, we are
guaranteed that a coagulation is going to occur in both copies. We
interpret this event as a success, and declare every other possibility
a failure. Hence if $G$ is a~geometric random variable with success
probability $n^{-1/4}/2$, and $(X_j)_{j=1}^\infty$ are i.i.d.
exponentials with mean $2n^{1/4}$, the total amount of time before
a~success occurs is dominated by $\sum_{j=1}^G X_j$.\vspace*{1pt}

If, however, one copy (say $\pi$) has only one unmatched part, then one first has to
break that component, which takes at most an exponential random
variable with rate $n^{-1/2}/4$. Note that the other copy must have had
at least~3 unmatched parts, so after breaking the big one, both copies
have now at least two unmatched copies and we are back to the preceding
case. It follows from this analysis
that in any case, $T_L$ is dominated by
\[
T_L \preceq Y+\sum_{j=1}^G X_j\vspace*{3pt}
\]
and so $\E(T_L) \le4 n^{1/2} + 4n^{1/2} = 8 n^{1/2}$.
Now, let
\[
\tau_L = T_L + T_{L-2} + \cdots+ T_4\vspace*{3pt}
\]
and let $T= \tau_{2n^{1/8}}$. Then $T$ is the time to get rid of all unmatched
parts. We obtain from the above
$\E(T) \le16 n^{5/8}$.
By Markov's inequality, it follows that $T< n^{5/8} \log n = \Delta$
with high probability. This concludes the proof of Theorem~\ref{T:mix}.
\end{pf*}

\section*{Acknowledgments}
We thank Hubert Lacoin, Remi Leblond and James Martin for a careful
reading of the first version of this manuscript and for constructive
comments and corrections. N. Berestycki is grateful to the Weizmann
Institute, Microsoft Research's Theory Group and the Technion for their
invitations in June 2008, August 2008 and April 2009, respectively.
Some of this research was carried out during these visits.


%
\printaddresses

\end{document}